# Chapter 1
# Lower Bounds for Identifying Subset Members with Subset Queries

E. Knill*

**Abstract**

An instance of a group testing problem is a set of objects $\mathcal{O}$ and an unknown subset $P$ of $\mathcal{O}$. The task is to determine $P$ by using queries of the type "does $P$ intersect $Q$", where $Q$ is a subset of $\mathcal{O}$. This problem occurs in areas such as fault detection, multiaccess communications, optimal search, blood testing and chromosome mapping. Consider the two stage algorithm for solving a group testing problem. In the first stage a predetermined set of queries are asked in parallel and in the second stage, $P$ is determined by testing individual objects. Let $n = |\mathcal{O}|$. Suppose that $P$ is generated by independently adding each $x \in \mathcal{O}$ to $P$ with probability $p/n$. Let $q_1$ ($q_2$) be the number of queries asked in the first (second) stage of this algorithm. We show that if $q_1 = o(\log(n)\log(n)/\log\log(n))$, then $\mathrm{Exp}(q_2) = n^{1-o(1)}$, while there exist algorithms with $q_1 = O(\log(n)\log(n)/\log\log(n))$ and $\mathrm{Exp}(q_2) = o(1)$. The proof involves a relaxation technique which can be used with arbitrary distributions. The best previously known bound is $q_1 + \mathrm{Exp}(q_2) = \Omega(p\log(n))$. For general group testing algorithms, our results imply that if the average number of queries over the course of $n^\gamma$ ($\gamma > 0$) independent experiments is $O(n^{1-\epsilon})$, then with high probability $\Omega(\log(n)\log(n)/\log\log(n))$ non-singleton subsets are queried. This settles a conjecture of Bill Bruno and David Torney and has important consequences for the use of group testing in screening DNA libraries and other applications where it is more cost effective to use non-adaptive algorithms and/or too expensive to prepare a subset $Q$ for its first test.

## 1 Introduction

An instance of a group testing problem is a set of objects $\mathcal{O}$ and an unknown subset $P$ of $\mathcal{O}$. The task is to determine $P$ by using queries of the type "does $P$ intersect $Q$", where $Q$ is an arbitrary subset (*pool*) of $\mathcal{O}$. An element $x \in \mathcal{O}$ is *positive* if $x \in P$, *negative* otherwise[1]. A pool is said to be positive if one of its objects is positive, negative otherwise. The determination of whether $P$ intersects $Q$ is called a *test* of $Q$. See Section 2 for a brief overview of the history and applications of group testing.

Algorithms for solving group testing problems can be classified by the degree to which they are adaptive. General adaptive algorithms can be modelled by an arbitrary binary decision tree, where each node corresponds to a pool $Q$, and which child to consider next is determined by the outcome of the test of $Q$. Completely non-adaptive algorithms (also called *one-stage* algorithms) are defined by a set of pools $\mathcal{Q}$, where all the pools in $\mathcal{Q}$ are tested in parallel and the set $P$ has to be determined from the outcome of the tests. A nearly non-adaptive algorithm that is of great interest for many screening problems is the *trivial two-stage* algorithm. Such an algorithm proceeds in two stages. In the first stage, the members of a fixed set of pools are tested in parallel and in the second stage only *individual* objects are tested. Which individual objects are tested may depend on the outcome of the first stage.

Research in the theory of group testing traditionally falls into two categories: probabilistic group testing and combinatorial group testing. In probabilistic group testing, a probabilistic model for the occurrence of positives is assumed, and the group testing procedure is optimized for minimum expected cost subject to constraints. In combinatorial group testing, it is assumed that the set of positives can be any member of a given family of sets $\mathcal{F}$, and the task is to find the algorithm which requires the minimum number of tests to uniquely determine a $P \in \mathcal{F}$ in the worst case. Combinatorial group testing is covered in detail in [6].

Here we resolve some questions about optimal algorithms for probabilistic group testing. Let the set of positive objects be distributed according to $\mathrm{Prob}(P)$, the probability that $P$ is the set of positives. The main contribution of this work is the development of a technique for obtaining lower bounds on the tradeoff between the number of pools $v = |\mathcal{Q}|$ and the expected number $\tilde{d}$ of queries that must be asked in the second stage. A fundamental result obtained by using this technique is as follows.

THEOREM 1.1. *Consider the group testing problem with $n$ objects and $v$ first-stage pools, where the set of*

[1] Positive elements are usually referred to as *defectives* in the literature. This choice of terminology is unfortunate, as in most applications of group testing today, no defect is implied.

positive objects $P$ is a uniformly randomly chosen $p$-tuple. Then the expected number $\tilde{p}$ of objects which occur only in positive pools is at least $(n/2^{v/p})(1 - \frac{4}{k}) - k$.

This is a restatement of Theorem 3.2 which is proved in Section 3.3. Note that if every $P \subseteq Q$, has even a very small chance of being the set of positive objects, then $\tilde{p} - p$ is a lower bound on the expected number of individual queries that must be asked in a trivial second stage to completely determine the positive objects.

We use Theorem 1.1 to obtain two results of practical interest. The first result applies to trivial two-stage group testing algorithms. Let the distribution of $P$ be determined by randomly and independently adding each object to $P$ with probability $p/n$ ($p$ constant). We say that $P$ is *Bernoulli* with parameter $p$.

THEOREM 1.2. *If $v \leq \beta^2 \ln(n) \log_2(n)/\ln\ln(n)$, then $\tilde{d} \geq n^{1-2\beta-o(1)}$. There are two-stage algorithms with $v = e^{1+\beta} \ln(n) \ln(n)/\ln\ln(n)(1 + o(1))$ and $\tilde{d} = O(n^{-\beta/2})$.*

Theorem 1.2 hints at a strong threshold behavior of $\tilde{d}$ given $v$.

The lower bounds obtained by our technique are substantially stronger than previously known ones. In the case where $\text{Exp}(|P|) = o(n)$, the best known bounds have been obtained by information theoretic arguments. The expected number of queries is at least $I(P) = \sum_{P \subseteq \mathcal{O}} \text{Prob}(P) \log_2(\text{Prob}(P))$, on average. For the case where $P$ is Bernoulli with parameter $p$, $I(P) = p \log_2(n)(1 + o(n))$. A notable feature of Theorem 1.2 is that the bounds are independent of $p$.

General algorithms for group testing can achieve the information theoretic bound to within a constant factor in the unit cost per query model. The simple bisection strategy finds $p$ positives in at most $2p\lceil \log_2(n) \rceil$ tests. If the distribution of $P$ is uniform over all $p$-tuples of $\mathcal{O}$ and $p = o(n)$, then the two-stage adaptive strategy also achieves this bound to within a constant factor ([4], Section 4). Theorem 1.2 shows that in general the two-stage strategy does not achieve this bound, even if $I(P) = O(\log_2(n))$. Intuitively, the reason for this is that the design of the pools must accommodate numbers of positives which are much larger than average, even though the probability of such an event is very small.

The second result derived from Theorem 1.1 applies to arbitrary adaptive algorithms. The lower bound of Theorem 1.2 implies that in $n^\gamma$ independent determinations of $P$ for the same set of objects, with high probability either the average number of individual object queries is $\Omega(n^{1-\epsilon})$ ($0 < \epsilon < 1$) or the number of non-singleton pools constructed is $\Omega(\log(n)\log(n)/\log\log(n))$. The details are given in Theorem 5.1.

## 2 Overview of group testing applications and significance of results

Group testing has been a much researched topic since the problem was formally published as a potential approach for economical blood testing [5]. In the blood testing problem, the task is to efficiently find the few samples which are positive for a disease such as syphilis by pooling samples and testing the pools. The basic idea is that if a pool is found to be negative, then all the samples which contributed to this pool can be excluded and do not have to be individually tested. This idea has since been used for quality control in product testing (when multiple items can be tested simultaneously) [20], searching files in punch card storage systems [15], efficient access of magnetic core memories [15], sequential screening of experimental variables [16], efficient algorithms for multiple-access systems and communication [17], and unique sequence screening of clone libraries [3, 4]. It has also been used in coding theory, optimal search, and the design of algorithms.

Adaptive and non-adaptive group testing algorithms seem to have been discovered and discussed independently. The first published paper on group testing [5] discusses a simple adaptive algorithm for probabilistic group testing. This paper gave rise to many studies of adaptive group testing algorithms both in probabilistic and combinatorial contexts (for example [20, 19, 18]). Non-adaptive group testing methods were discovered somewhat later in the context of efficient searching of punch card files and accessing magnetic core memories (see [15] and the references therein). The first methods used for these applications were randomized and and used primarily in a probabilistic context. Kautz and Singleton [15] first used combinatorial methods from coding theory to obtain combinatorial non-adaptive group testing algorithms. Kautz and Singleton's work was continued by Dyachov and others in the Soviet Union [7, 8, 9, 10, 11] and elsewhere [13, 14, 21, 17, 12, 3, 1, 4].

The results of this paper have immediate implications for applications of group testing where algorithms with more than two stages are undesirable and/or pools can be reused but the initial construction of non-singleton pools is expensive compared to testing. In most previous studies of group testing, the cost model assumes unit cost per query and does not consider pool construction independently from pool testing or testing of individual objects. This is appropriate for problems where pools cannot be reused or are cheap to build. It does not address the actual costs encountered in screening clone or protein libraries, which are possibly the most active current applications of group testing.

Nearly every laboratory involved in the mapping of chromosomes uses group testing for library screening. In this application, the chromosome or genome of interest (basically a sequence of DNA) is randomly cut into many overlapping pieces of similar sizes. These pieces are replicated in clones and stored in large genomic libraries (1000 to 75000 clones, there has been discussion of libraries with up to $10^6$ clones). The first task is to determine the arrangement of these clones on the original sequence of DNA. One of the methods for doing this is to obtain for each of a large set of unique sites the set of clones which contain this site. Once this is done, the sites and clones are ordered and localized[2]. One would then like to use the clone libraries to find genes and other features of interest. This also requires screening the clones.

There are two features which distinguish the library screening problem from many other applications of group testing. The first is that the same library will be tested for many sites. In fact, the number of sites that must be tested before the sequence can be reliably reconstructed is usually of the order of the number of clones in the library. This is why Theorem 5.1 is relevant. The second feature is that pool construction is very costly. It is generally feasible to construct a number of pools (much fewer than the number of clones) initially by exploiting parallelism, but adaptive construction of pools with many clones during the testing procedure is discouraged. The technicians who implement the pooling strategies generally dislike even the 3-stage strategies that are often used. Thus the most commonly used strategies for pooling libraries of clones rely on a fixed but reasonably small set of non-singleton pools. The pools are either tested all at once or in a small number of stages (usually at most 2) where the previous stage determines which pools to test in the next stage. The potential positives are then inferred and confirmed by testing of individual clones. In most biological applications each positive clone must be confirmed even if the pool results unambiguously indicate that it is positive. This is to improve the confidence in the results, given that in practice the tests are prone to errors.

The first formal study of how to best screen libraries of clones is due to Barillot et al. [3]. They showed that fairly efficient trivial two stage strategies can be obtained by using simple geometric constructions for the pools. It has since been realized [2] that these constructions correspond to simple error correcting codes and were already discussed in [15]. In [4] randomized constructions are shown to be very useful for library screening. These types of randomized constructions were originally used for the file searching application [15] and were thoroughly analyzed by Dyachov [8] in a combinatorial group testing context.

For screening libraries of clones, the distribution of $P$, the set of positive clones, is well approximated by the Bernoulli distribution, at least for the first $O(n)$ screenings. Recently computational studies of W.J. Bruno (unpublished) showed that for the Bernoulli distribution the number of pools required for successfully obtaining the positive clones with few individual tests appeared to grow substantially faster for randomized constructions than the information theoretic lower bound. He conjectured that such a growth was necessary for the Bernoulli distribution. The results of this paper imply that this conjecture is true.

The results given here still hold provided that the distribution of positives has a sufficiently large probability of $\Omega(\log(n)/\log\log(n))$ positives (see the proof of Theorem 3.3). This means that the effects of the bound will be observable provided that the distribution is close to a Bernoulli distribution, as is the case for the initial screenings of a library of clones. This is due to the random nature of the construction of clone libraries. In a fixed library, the information obtained from previous screenings will eventually constrain the possible results. However, unless the ordering of the tested sites is known in advance, the number of screenings required to observe this is $\Omega(n)$. Other group testing applications such as blood testing or affinity testing for proteins may show fewer dependencies between screenings. It is therefore likely that the tradeoffs of Theorems 1.2 and 5.1 can be observed in practice. Provided that it is indeed desirable to completely determine the set of positive objects, the number of non-singleton pools that must be constructed is substantially larger than required by the information theoretic bound.

The remainder of this work is organized as follows: Section 3 describes a general technique for obtaining lower bounds for non-adaptive group testing strategies. The lower bound of Theorem 1.2 follows immediately from Theorem 3.4 which is proved at the end of this section. In Section 4 we use the probabilistic method to obtain the upper bound of Theorem 1.2. This bound is a consequence of Theorem 4.1. In Section 5 it is shown how Theorem 1.2 can be used to obtain bounds on the minimum number of non-singleton pools that are constructed by any algorithm. Some open problems and directions for future work are given in Section 6.

---

[2]How to obtain the reconstruction from this (usually imperfect) data is itself the topic of intensive research in approximation algorithms.

## 3 Lower bound methods

The technique for obtaining lower bounds on the trade-off between the number of pools and the number of second stage queries in two-stage algorithms relies on a combinatorial relaxation technique which transforms the problem to a linear program. The initial linear program is simplified for symmetric distributions so that in principle it could be solved exactly. We estimate its optimum value by going to the dual and applying a simple greedy method to obtain good bounds.

**3.1 Relaxation to a linear program.** Consider the general problem of constructing informative pools. Let $\mathcal{Q}$ be a set of $v$ pools. Write $S = 2^v$. Any function $g : \mathcal{O} \to 2^{\mathcal{Q}}$ determines a way of pooling the objects by adding each object $x$ to the pools in $g(x)$. If $P$ is the set of positive objects, then the set of positive pools is given by

$$\bigcup_{x \in P} g(x) = \{y \mid \exists x \in P \text{ such that } y \in g(x)\}.$$

Let $\text{Prob}(P)$ be the probability that $P \subseteq \mathcal{O}$ is the set of positive objects.

Suppose all the pools of $\mathcal{Q}$ are tested and the set of positive pools is $\mathcal{Q}_p$. An object $x$ can be positive only if $g(x) \subseteq \mathcal{Q}_p$. Such objects are called *candidate positives*. A candidate positive object $x$ with $x \notin P$ is called *unresolved negative*. Let $\tilde{d}$ be the expected number of unresolved negative objects. In many cases of interest, at least the unresolved negative objects must be examined by any second stage individual testing method. This occurs in particular if $\text{Prob}(P) > 0$ for each $P$, which is satisfied by the Bernoulli distribution. In general, $\tilde{d}$ is a good estimate of the number of second stage tests whenever the distribution is sufficiently rich. Note that in practice, it is often the case that all candidate positive objects of interest are confirmed negative or positive to improve the confidence (for some alternative approaches, see [4]).

Let $\bar{P}$ be the set of candidate positive objects,

$$\bar{P} = \{y \in \mathcal{O} \mid g(y) \subseteq \bigcup_{x \in P} g(x)\}.$$

The operation $P \to \bar{P}$ is a closure operation which frequently occurs in the study of union-closed families of sets and lattices. The expected number of unresolved negative objects $\tilde{d}$ is computed as

$$(3.1) \qquad \tilde{d} = \sum_{P \subseteq \mathcal{O}} \text{Prob}(P) |\bar{P} \setminus P|.$$

The goal is a lower bound on $\tilde{d}$ for the given distribution and number of pools by minimizing $\tilde{d}$ over all functions $g$. As it stands, this optimization problem is difficult. We can however relax the problem to linear programming by allowing physically impossible pools.

To obtain the desired linear program $L(n, S)$, we shift the problem to $\mathcal{O}$. The choice of pools is replaced first by a choice of a suitable closure operation $P \to \bar{P}$. A subset of $\mathcal{O}$ is *closed* if it is given by $\bar{P}$ for some $P \subseteq \mathcal{O}$. The closed subsets of $\mathcal{O}$ form an intersection-closed family of sets which contains $\mathcal{O}$. Two important observations are: (1) The number of closed subsets is at most $S$, since every subset $U$ of the pools determines the closed set $\{x \mid g(x) \subseteq U\}$ and all closed sets can be obtained like this. (2) For each $P$, $\bar{P}$ is the unique minimal closed subset which includes $P$.

A weak fractional version of the closed sets can be described by a set of variables $w_V$ and $w_{U,V}$ for $U \subseteq V \subseteq \mathcal{O}$ subject to the constraints

$$(3.2) \qquad w_V \geq 0, \quad w_{U,V} \geq 0,$$

$$(3.3) \qquad \sum_V w_V \leq S,$$

$$(3.4) \qquad \text{for each } U: \sum_{V:U \subseteq V} w_{U,V} \geq 1,$$

$$(3.5) \qquad \text{for each } U \subseteq V: \quad w_{U,V} \leq w_V.$$

The closed sets correspond to the variables $w_V$ and the cardinality constraint is enforced by inequality (3.3). How much of each "closed" subset corresponds to the unique minimal one which includes a given $V$ is now described by the variables $w_{U,V}$. Inequalities (3.5) ensure that the amount of $V$ which includes $U$ does not exceed the degree to which $V$ is "closed". Inequalities (3.4) ensure that each $U$ is included in a total of at least one "closed" subset.

Formally, a feasible solution of $L(n, S)$ can be obtained from a $g : \mathcal{O} \to \mathcal{P}$ by defining

$$\begin{aligned} w_U &= [\bar{U} = U], \\ w_{U,V} &= [V = \bar{U}], \end{aligned}$$

where for any logical expression $\phi$, $[\phi] = 1$ if $\phi$ is true, and $[\phi] = 0$ otherwise.

A lower bound on $\tilde{d}$ is now determined by the minimum value $\text{val}(L(n, S), \text{Prob})$ of

$$(3.6) \qquad \sum_U \text{Prob}(U) \sum_{V:U \subseteq V} w_{U,V} |V \setminus U|.$$

Given that the number of variables is $2^n + 3^n$, the problem of evaluating $\text{val}(L(n, S), \text{Prob})$ is impractical in general. However in many cases of interest, the probability distribution is symmetric, that is, $\text{Prob}(U)$ depends only on $|U|$. Since $L(n, S)$ itself is symmetric, the number of variables can be substantially reduced by assuming that $w_V$ and $w_{U,V}$ depend only on the cardinalities of $U$ and $V$. In particular, in the symmetric case, $L(n, S)$ is equivalent to $L'(n, S)$ with variables $w_j$ and $w_{i,j}$ for $0 \leq i \leq j \leq n$ and constraints

$$(3.7) \qquad w_j \geq 0, \quad w_{i,j} \geq 0,$$

$$(3.8) \qquad \sum_j \binom{n}{j} w_j \leq S,$$

$$(3.9) \quad \text{for each } i: \quad \sum_{j:j \geq i} \binom{n-i}{j-i} w_{i,j} \geq 1,$$

$$(3.10) \text{ for each } i \leq j: \quad w_{i,j} \leq w_j.$$

The quantity to be minimized is

$$(3.11) \qquad \sum_i p(i) \sum_{j:j \geq i} \binom{n-i}{j-i} w_{i,j}(j-i),$$

where $p(i) = \sum_{|U|=i} \text{Prob}(U)$.

To find useful lower bounds on $\text{val}(L'(n, S), \text{Prob})$, we can use linear programming duality. The dual program $L^*(n, S)$ has variables $v$ (for inequality (3.8)), $v_i$ ($0 \leq i \leq n$, for inequalities (3.9)) and $v_{i,j}$ ($0 \leq i \leq j \leq n$, for inequalities (3.10)). The constraints are

$$(3.12) \qquad v \geq 0, \quad v_i \geq 0, \quad v_{i,j} \geq 0,$$

$$(3.13) \qquad \text{for each } j: \sum_{i:i \leq j} v_{i,j} - \binom{n}{j} v \leq 0,$$

(3.14) for each $i \leq j$:
$$(3.15) \quad \binom{n-i}{j-i} v_i - v_{i,j} \leq p(i) \binom{n-i}{j-i}(j-i).$$

The value $\text{val}(L^*(n, S), \text{Prob})$ is given by the maximum value of
$$(3.16) \qquad l^* = \sum_i v_i - Sv.$$

**3.2 A simple feasible solution for $L^*$.**
By the duality theorem of linear programming, $\text{val}(L'(n, S), \text{Prob}) = \text{val}(L^*(n, S), \text{Prob})$ and any feasible solution yields a lower bound. We use what amounts to a greedy method to find a reasonably good solution. The idea is to increase each $v_i$ and adjust the other variables as long as $l^*$ increases as a result. If one of the inequalities (3.14) is violated, compensate by increasing the values of the appropriate $v_{i,j}$. To satisfy inequalities (3.13) may require increasing $v$. The change of $v_i$ is successful at increasing $l^*$ if the adjustment of $v$ is sufficiently small.

For each $i$, let

$$s(i) = \min\{s \mid (s)_i > (n)_i/S\},$$

where $(s)_i = s(s-1)...(s-i+1)$ is the $i$'th falling factorial of $s$.

THEOREM 3.1.
$$\text{val}(L^*(n, S), \text{Prob}) \geq$$
$$\sum_i p(i)(s(i) - i) -$$
$$S \max_j \sum_{i:i \leq j < s(i)} \frac{(j)_i}{(n)_i} p(i)(s(i) - j).$$

*Proof.* By following a greedy strategy of finding a good feasible solution, one obtains

$$v_i = p(i)(s(i) - i).$$

To extend this to a feasible solution of $L^*(n, S)$, let

$$v_{i,j} = \binom{n-i}{j-i} p(i)(s(i) - j) \ [i \leq j < s(i)].$$

This ensures that inequalites (3.14) are satisfied. To satisfy inequalites (3.13) requires that for each $j$,

$$\binom{n}{j} v \geq \sum_{i:i \leq j} v_{i,j} = \sum_{i:i \leq j} [j < s(i)] \binom{n-i}{j-i} p(i)(s(i) - j).$$

Thus we can let

$$v = \max_j \sum_{i:i \leq j < s(i)} \frac{(j)_i}{(n)_i} p(i)(s(i) - j).$$

For $l^*$ we now get

$$l^* = \sum_i p(i)(s(i) - i) - S \max_j \sum_{i:i \leq j < s(i)} \frac{(j)_i}{(n)_i} p(i)(s(i) - j).$$
(3.17)

The right-hand side of the inequality in Theorem 3.1 can be reasonably estimated for any symmetric distribution and gives a lower bound for the optimum value of $\tilde{d}$ given $S$.

**3.3 Evaluation of the bound for the Bernoulli distribution.** For the distributions of interest here, we can obtain a lower bound on $\tilde{d}$ by conditioning on the case where the number of positives is fixed. Thus we estimate $\tilde{d}$ by considering the distribution $p_k(i) = [i = k]$. Let $b_k(n, S) = \text{val}(L^*(n, S), p_k)$. The minimum expected number of unresolved negatives can then be estimated by $p(k) b_k(n, S)$.

LEMMA 3.1. *For all $n$, $S$ and $i$, $s(i) > \frac{n}{S^{1/i}}$.*

*Proof.* We have $(s(i))_i > (n)_i/S$. Since $s(i) \leq n$, $(s(i))_i/(n)_i \leq (s(i)/n)^i$ and the result follows.

THEOREM 3.2.
$$b_k(n,S) \geq s(k)(1-4/k) - k.$$

*Proof.* We can assume without loss of generality that $s(k) \geq k$ and $k \geq 2$. We have
$$b_k(n,S) \geq s(k) - k - S \max_{j:k\leq j<s(k)} \frac{(j)_k}{(n)_k}(s(k)-j).$$

Observe that by the definition of $s(k)$, $S/(n)_k \leq 1/(s(k)-1)_k$. Hence
$$b_k(n,S) \geq s(k) - k - \frac{1}{(s(k)-1)_k} \max_{j:k\leq j<s(k)} (j)_k(s(k)-j).$$

Write $s = s(k)$ and let $f(j) = (j)_k(s-j)$. The maximum of $f(j)$ occurs at a $j$ such that $f(j)/f(j-1) \geq 1$ and $f(j+1)/f(j) < 1$. Since $f(j)/f(j-1) = (s-j)j/((s-j+1)(j-k))$, the following are equivalent:
$$\begin{aligned}
f(j)/f(j-1) &\geq 1 \\
sj - j^2 &\geq -(s+1)k + (k+s+1)j - j^2 \\
(s+1)k &\geq (k+1)j \\
(s+1)k/(k+1) &\geq j.
\end{aligned}$$

This implies that the maximum occurs at $j_m = \min(s-1, \lfloor(s+1)k/(k+1)\rfloor)$.

$$\begin{aligned}
f(j_m)/(s-1)_k &\leq (\lfloor(s+1)(k/(k+1))\rfloor/(s-1))^k (s-j_m) \\
&\leq ((s+1)/(s-1))^k 1/(1+1/k)^k(s+1)/(k+1) \\
&\leq ((1+1/k)/(1-1/k))^k 1/(1+1/k)^k(s/k) \\
&\leq 4(s/k).
\end{aligned}$$

where we used $s \geq k \geq 2$ and the fact that $(1-1/k)^k$ is increasing in $k$. The proof of the theorem is completed by substituting this inequality for the last summand in the lower bound for $b_k(n,S)$.

THEOREM 3.3. *Let $\alpha > \beta > 0$, $k = \alpha f(n)(1+o(1))$ and $v \leq \beta \log_2(n) f(n)$, with $f(n) = o(n^\epsilon)$ for some $0 < \epsilon < 1 - \beta/\alpha$ and $1 = o(f(n))$. Then*
$$b_k(n,S) \geq n^{1-\beta/\alpha-o(1)}.$$

*Proof.* Applying the previous results gives
$$\begin{aligned}
b_k(n,S) &\geq s(k)(1-o(1)) - k \\
&\geq n/S^{1/k}(1-o(1)) - k \\
&= n^{1-\beta/\alpha-o(1)}.
\end{aligned}$$

For the rest of this section, let $p(i)$ be determined by the Bernoulli distribution with parameter $p$, so that $p(i) = \binom{n}{i}(\frac{p}{n})^i(1-\frac{p}{n})^{n-i}$. The following theorem implies the first half of Theorem 1.2.

THEOREM 3.4. *Let $0 < \beta < 1/2$. For $v \leq \beta^2 \ln(n) \log_2(n)/\ln\ln(n)$,*
$$\tilde{d} \geq n^{1-2\beta-o(1)}.$$

*Proof.* We can estimate $\tilde{d}$ by $p(k)b_k(n,S)$ with $k = \alpha \ln(n)/\ln\ln(n)(1+o(1))$. The probability $p(k)$ is bounded as follows:
$$\begin{aligned}
p(k) &= \binom{n}{k}\left(\frac{p}{n}\right)^k\left(1-\frac{p}{n}\right)^{n-k} \\
&\geq \frac{1}{k^k}e^{-O(k)} \\
&= e^{-k\ln(k)(1+o(1))} \\
&= n^{-\alpha(1+o(1))}.
\end{aligned}$$

Using Theorem 3.3, this gives
$$\tilde{d} \geq p(k)b_k(n,S) \geq n^{-\alpha+1-\beta^2/\alpha-o(1)}.$$

Let $\alpha = \beta$. Then $\tilde{d} \geq n^{1-2\beta-o(1)}$.

It is clear that a result such as Theorem 3.4 holds for any distribution $p$ which satisfies $p(k) \geq n^{-\alpha-o(1)}$ for suitable $k$. In fact, any $k$ for which this holds implies a tradeoff similar to that in Theorem 1.2.

## 4 Probabilistic construction of pools

The probabilistic method can be used to show that the pools can be chosen in a nearly optimal fashion. Let $\mathcal{P}$ be the set of $v$ pools. In this section we specify the relationships between the objects and the pools by an $n \times v$ incidence matrix $I$, where $I_{x,y} = 1$ if object $x$ occurs in pool $y$ and $I_{x,y} = 0$ otherwise. Consider $I$ as a random variable with distribution determined by randomly and independently setting each $I_{x,y}$ to 1 with probability $q$. This and related models are frequently used for probabilistic constructions involving general families of sets and have been applied to one-stage non-adaptive group testing algorithms [7, 8, 9, 12].

Let $\tilde{d}(I)$ be the expected number of unresolved negatives for the pools constructed according to $I$ if the probability of exactly $i$ positives is $p(i)$. The following calculations can be done for more general distributions of the positive objects to obtain similar results.

LEMMA 4.1.
$$\tilde{d}(I) = \sum_{i=0}^n p(i)(n-i)(1-q(1-q)^i)^v.$$

*Proof.* Assume that there are $i$ positive objects $x_1, \ldots, x_i$, and $x$ is a negative object. For $x$ to be unresolved negative requires that for every pool $y$ the following is not the case: $x$ is in $y$ (probability $q$) and for each $x_i$, $x_i$ is not in $y$ (probability $(1-q)^i$). Thus the probability that $x$ is unresolved negative is given by $(1-q(1-q)^i)^v$. Given that there are $i$ positive objects, there are $n-i$ potential unresolved negative objects. Using linearity of expectations yields the sum in the lemma.

Consider now the case where the distribution $p(i)$ is determined by each object's being independently positive with probability $p/n$.

LEMMA 4.2.
$$\tilde{d}(I) \leq n \sum_{i=0}^{\infty} \frac{p^i}{i!} e^{-q(1-q)^i v}.$$

*Proof.*
$$\begin{aligned}
\tilde{d}(I) &= \sum_{i=0}^{n} p(i)(n-i)(1-q(1-q)^i)^v \\
&\leq n \sum_{i=0}^{n} \binom{n}{i} \left(\frac{p}{n}\right)^i \left(1-\frac{p}{n}\right)^{n-i} (1-q(1-q)^i)^v \\
&\leq n \sum_{i=0}^{n} \frac{p^i}{i!} (1-q(1-q)^i)^v \\
&\leq n \sum_{i=0}^{\infty} \frac{p^i}{i!} e^{-q(1-q)^i v}.
\end{aligned}$$

THEOREM 4.1.
Let $\epsilon > 0$, $v = e^{1+\epsilon} \ln(n) \ln(n)/\ln \ln(n)(1 + o(1))$ and $q = \ln \ln(n)/\ln(n)$. Then
$$\tilde{d}(I) = O(n^{-\epsilon/2}).$$

*Proof.* Let $k = (1+\epsilon/2+\delta) \ln(n)/\ln \ln(n)$ with $\delta > 0$ constant but sufficiently small as required by the calculations below. Divide the sum of Lemma 4.2 into two parts,
$$\begin{aligned}
S_1 &= n \sum_{i: 0 \leq i \leq k} \frac{p^i}{i!} e^{-q(1-q)^i v}, \\
S_2 &= n \sum_{i: i > k} \frac{p^i}{i!} e^{-q(1-q)^i v}.
\end{aligned}$$

Estimate $S_2$ by
$$\begin{aligned}
S_2 &\leq n \sum_{i: i > k} \frac{p^i}{i!} \\
&\leq n \frac{p^k}{k!} e^p
\end{aligned}$$

$$\begin{aligned}
&= n e^{-k \ln(k)(1-o(1))} \\
&= n^{(-\epsilon/2-\delta)(1-o(1))} \\
&= O(n^{-\epsilon/2}).
\end{aligned}$$

Bound $S_1$ as follows:
$$\begin{aligned}
S_1 &\leq n e^{-q(1-q)^k v} e^p \\
&= n e^{-\exp(1+\epsilon) \exp(-1-\epsilon/2-\delta)(1-o(1)) \ln(n)} \\
&= n^{1-\exp(\epsilon/2-\delta)(1-o(1))} \\
&= O(n^{-\epsilon/2}),
\end{aligned}$$

where we have used the fact that
$$\begin{aligned}
(1-q)^k &= \left(1 - \frac{1+\epsilon/2+\delta}{k}\right)^k \\
&= e^{-1-\epsilon/2-\delta}(1 - O(1/k)).
\end{aligned}$$

We now have $\tilde{d}(I) \leq S_1 + S_2 = O(n^{-\epsilon/2})$.

## 5 Implications for general algorithms

To see what Theorem 1.2 implies for general algorithms applied to independent instances of the Bernoulli distribution requires estimating the probability that $\tilde{d} > n^\epsilon$.

LEMMA 5.1. *Under the assumptions of Theorem 3.4, $\text{Prob}(\tilde{d} \geq n^{1-2\beta-\epsilon}) \geq n^{-2\beta-o(1)}$ for any fixed $\epsilon > 0$.*

*Proof.* Since $\tilde{d} \leq n$ we can use the Chebyshev inequality as follows:
$$\begin{aligned}
\text{Prob}(n - \tilde{d} \geq n - n^{1-2\beta-\epsilon}) &\leq \frac{\text{Exp}(n-\tilde{d})}{n - n^{1-2\beta-\epsilon}} \\
&\leq \frac{n - n^{1-2\beta-o(1)}}{n - n^{1-2\beta-\epsilon}}.
\end{aligned}$$

Hence
$$\begin{aligned}
\text{Prob}(\tilde{d} \geq n^{1-2\beta-\epsilon}) &\geq 1 - \frac{n - n^{1-2\beta-o(1)}}{n - n^{1-2\beta-\epsilon}} \\
&= \frac{n^{1-2\beta-o(1)} - n^{1-2\beta-\epsilon}}{n - n^{1-2\beta-\epsilon}} \\
&= n^{-2\beta-o(1)}.
\end{aligned}$$

THEOREM 5.1. *Let $0 < \gamma$, $0 < \beta < 1/4$ and $2\beta < \gamma$. Suppose that a group testing algorithm is applied to $n^\gamma$ independent Bernoulli instances of $P$. Then with probability at least $1 - e^{-n^{\gamma-2\beta-o(1)}}$, either the average number of individual tests is greater than $n^{1-4\beta-o(1)}$ or more than $\beta^2 \ln(n) \log_2(n)/\ln \ln(n)$ non-singleton pools are constructed.*

*Proof.* Let $\epsilon > 0$ be arbitrarily small. Let $E$ be the event which consists of constructing a new non-singleton pool, or individually testing at least $n^{1-2\beta-\epsilon/2}$ objects

, or already having $> \beta^2 \ln(n) \ln(n)/\ln\ln(n)$ pools. By Theorem 1.2 and Lemma 5.1, for each instance of $P$ the probability of $E$ is at least $n^{-2\beta - o(1)}$. Let $n(E)$ be the number of events $E$ that occur in $n^\gamma$ instances. We would like to relate the distribution of $n(E)$ to a binomial distribution and apply tail estimates for the binomial distribution. That this can be done follows from the next lemma.

LEMMA 5.2. *Let $E \subseteq \{1, \ldots, k\}$ be a random variable which satisfies that for each $U$, $Prob(i \in E \mid E \cap \{1, \ldots, i-1\} = U) \geq q$. Let $B \subseteq \{1, \ldots, k\}$ be the Bernoulli random variable with $Prob(i \in B) = q$. Let $\mathcal{F}$ be an upward closed family of subsets of $\{1, \ldots, k\}$ (i.e. $U \in \mathcal{F}$ and $V \supset U$ implies $V \in \mathcal{F}$). Then $Prob(E \in \mathcal{F}) \geq Prob(B \in \mathcal{F})$.*

*Proof.* We use induction on $k$. The lemma holds trivially for $k = 1$. Let $\mathcal{F}_0 = \{U \mid U \in \mathcal{F} \wedge U \subseteq \{2, \ldots, k\}\}$ and $\mathcal{F}_1 = \{U \cap \{2, \ldots, k\} \mid U \in \mathcal{F}\}$. By upward closure, $\mathcal{F}_0 \subseteq \mathcal{F}_1$.

$$\begin{aligned}
&\text{Prob}(E \in \mathcal{F}) \\
&= \text{Prob}(1 \in E \wedge E \cap \{2, \ldots, k\} \in \mathcal{F}_1) \\
&\quad + \text{Prob}(1 \notin E \wedge E \in \mathcal{F}_0) \\
&= \text{Prob}(1 \in E)\text{Prob}(E \cap \{2, \ldots, k\} \in \mathcal{F}_1 \mid 1 \in E) \\
&\quad + \text{Prob}(1 \notin E)\text{Prob}(E \in \mathcal{F}_0 \mid 1 \notin E) \\
&\geq \text{Prob}(1 \in E)\text{Prob}(B \cap \{2, \ldots, k\} \in \mathcal{F}_1) \\
&\quad + \text{Prob}(1 \notin E)\text{Prob}(B \cap \{2, \ldots, k\} \in \mathcal{F}_0),
\end{aligned}$$

where the last step used the induction hypothesis twice for $\{2, \ldots, k\}$ with $E' = E \cap \{2, \ldots, k\}$ conditioned on $1 \in E$, and $E'' = E$ conditioned on $1 \notin E$. The result follows by using the inequalities $\text{Prob}(B \cap \{2, \ldots, k\} \in \mathcal{F}_1) \geq \text{Prob}(B \cap \{2, \ldots, k\} \in \mathcal{F}_0)$ and $\text{Prob}(1 \in E) \geq q$.

Using Lemma 5.2 with $\mathcal{F} = \{U \mid |U| > n^{\gamma - 2\beta - o(1)}\}$ allows us to apply the usual tail estimates on the binomial distribution to obtain

$$\text{Prob}(n(E) \leq n^{\gamma - 2\beta - o(1)}) \leq e^{-n^{\gamma - 2\beta - o(1)}},$$

where constants have been absorbed into $n^{-o(1)}$. Suppose that less than $n^{1 - 4\beta - \epsilon}$ individual tests are performed on average and less than $\beta^2 \ln(n) \ln(n)/\ln\ln(n)$ non-singleton pools are constructed. Then

$$n(E) \leq \beta^2 \ln(n) \ln(n)/\ln\ln(n) + n^{\gamma - 2\beta - \epsilon/2}$$

for the maximum number of times additional non-singleton pools are constructed or at least $n^{1 - 2\beta - \epsilon/2}$ objects are tested. It follows that $1 - e^{-n^{\gamma - 2\beta - o(1)}}$ is a lower bound on the probability of the event in the theorem.

## 6 Some problems

Some of the interesting questions that are raised by this work include:

PROBLEM 6.1. Determine the precise nature of the threshold behavior of the tradeoff between $v$ and $\tilde{d}$ for two-stage algorithms.

PROBLEM 6.2. Given the information $I(P)$ of the distribution of positives, what is the maximum gap between the the information theoretic lower bound and the number of pools and individual queries required by a two-stage algorithm?

PROBLEM 6.3. Consider an arbitrary symmetric distribution of positives. Is it true that up to a multiplicative constant the optimal two-stage algorithm is obtained by the probabilistic construction of Section 4?

PROBLEM 6.4. Can similar lower bounds be proved for approximation algorithms, that is algorithms which either determine $P$ with high probability, or find at least $\min(p, |P|)$ positives with high probability?

Problem 6.4 is suggested by work described in [4]. Note that if we are allowed to fail to determine $P$ in $n^{-o(1)}$ instances, then the tradeoffs in Theorem 1.2 do not apply.